%% file: 101204-reg-powers-ulrich.tex
\documentclass[twoside,10pt]{article}
\usepackage{amsmath,amsthm, amssymb,amsxtra}
\usepackage{graphicx}
\usepackage{diagrams}

\input preamble.tex

\newarrow{Into} C--->
\date{November 14, 2010}
\title{Stabilization of the Regularity of Powers of An Ideal
\footnote{The authors are grateful for the support of the National Science Foundation during the preparation of this work.}} 
\author{David Eisenbud and Bernd Ulrich}

\begin{document}

\maketitle

\begin{abstract}
When $M$ is a finitely generated graded module over a standard graded algebra $S$ and $I$ is an ideal of $S$, it is known from work of
Cutkosky, Herzog, Kodiyalam, R\"omer, Trung and Wang that the Castelnuovo-Mumford regularity of $I^mM$ has the form $dm+e$  when $m\gg0$. We 
give an explicit bound on the $m$ for which this is true, under the hypotheses that $I$ is generated in a single degree and $M/IM$ has 
finite length, and we explore the phenomena that occur when these hypotheses are not satisfied. Finally, we prove a regularity bound for 
a reduced, equidimensional projective scheme of codimension 2 that is similar to the bound in the Eisenbud-Goto conjecture [1984], under the 
additional hypotheses that the scheme lies on a quadric and has nice singularities.

\end{abstract}

\section*{Introduction}

Let $S$ be a standard graded algebra over a field $k$---that is, an algebra generated by finitely many forms of degree one---and 
let $M$ be a finitely generated graded $S$-module. If $H$ is an artinian $S$-module we set $\reg H = \max\{d\mid  H_d \neq 0\}$
 and we write $\reg M$ for the Castelnuovo Mumford regularity
$$
\reg M = \reg_{S_+}M := \max \{ \reg H^i_{S_+} + i\}.
$$ 
Combining results of Cutkosky-Herzog-Trung [1999], Kodiyalam [2000], R\"omer [2001] and Trung-Wang [2005], we have:
\begin{theorem}\label{CHTK}
There exist integers $m_0 = m_0(I,M)$, $d = d(I,M)$ and $e = e(I,M)$ such that for all $m\geq m_0$, 
$$
\reg I^m M = dm+e  \, .
$$
Furthermore, $d$ is the asymptotic generator degree of $I$ on $M$, i.e., the minimal number such that if $J\subset I$ is 
the ideal generated by the elements of $I$ of degree $\leq d$, then
$I+\ann M$ is integral over $J+\ann M$.
\end{theorem}

This beautiful result begs for an answer to several questions: What is the significance of the number $e$? What is a reasonable 
bound $m_0$? What is  the nature of the function $m\mapsto \reg I^mM$ for $m<m_0 \ \dots$? In general very little is known. But 
the result of the first section of this paper gives a value for $m_0$ in case
$$
(*)\quad \hbox{$I$ is generated in a single degree and $M/IM$ has finite length.}
$$
Here is a summary of our knowledge in this case. Under the hypothesis $(*)$ one has:
\begin{itemize}
\item The number $d$ in Theorem \ref{CHTK} is equal to the common degree of the generators of $I$.
\item The differences $e_m := \reg I^m -dm$ form a weakly decreasing sequence of non-negative integers.
\item The asymptotic value $e$ of the $e_m$ can be identified with the regularity of the restriction of
the sheaf associated to $M$ to the fibers of the morphism defined by $I$.
\item The numbers $e_m$ are equal to the assymptotic value $e$ for all $m\geq m_0$, where $m_0$ is the \emph{$(0,1)$-regularity} (defined below)
of the Rees algebra $\cR(I)$.
 \end{itemize}
The first item in this list is immediate from the definitions. The next two are proved in Eisenbud-Harris [2008]. The last is the subject 
of the first section of this paper, where we also derive a sharper but more technical bound that is often optimal. We note that a different (somewhat larger) value for $m_0$ was proposed in Cutkosky-Herzog-Trung [1999], but the proof given was incomplete, as the authors of that paper have pointed out. Marc Chardin has informed us that, after seeing our work, he was able to extend our bound on the power $m_0$. He uses a spectral sequence argument to treat the case of an ideal $I$ such that $M/IM$ has finite length, without assuming that $I$ is generated in a single degree.

In connection with the second item of the list, we observed in many cases that the sequence of first differences of the $e_{m} - e_{m+1}$ is also weakly 
decreasing. Is this always the case, under the assumption of $(*)$ ?

A key definition in this development is the $(0,1)$ (Castelnuovo-Mumford) regularity of the Rees module $\cR(I.M)$.
To define it, we recall that the Rees ring of $I$ is
$$
\cR(I) : = \oplus_{n\geq 0} I^n  \cong \oplus_{n\geq 0} I^nt^n =  S[It] \subset S[t] \, .
$$
This ring is an epimorphic image of the polynomial ring $T:=S[y_0, \ldots, y_r]$ via the map of $S$-algebras sending
the $y_i$ to $t$ times the homogeneous minimal generators of $I$. In fact, this becomes a map of bigraded 
$k$-algebras if we set
$\deg x_i = (1,0)$ and $\deg y_i = (0,1)$ (note that this is only possible because the generator degrees of $I$ are assumed to be equal).
Next, if $M$ is a finitely generated graded $S$-module, we define
$$
\cR (I, M) = S[It]M \subset M\otimes_SS[t] \, ,
$$ 
which is a finitely generated bigraded module over
$\cR(I)$ and hence over $T$. Thus we consider a bigraded minimal free resolution
$$
\cdots  \ \  F_1 \longrightarrow  F_0 \longrightarrow  \cR(I,M) \to 0
$$ 
of $\cR(I,M)$ as $T$-module, and we define
$\reg_{(0,1)} \cR(I,M)$ to be the maximum integer $j$ such that $F_i$ has a free summand of the form
$T(-a, -i-j)$ for some $i$ and $a$. As with the usual Castelnuovo-Mumford regularity, there is also a 
definition in terms of local cohomology, which we will use freely; see R\"omer [2001] for a detailed treatment.

In the second section of this paper, we turn to the question of what happens if we weaken the hypothesis $(*)$ to allow ideals 
that are not necessarily generated in a single degree. We found it surprisingly hard to give formulas for the numbers 
$e_m(I,M) := \reg I^mM - d(I,M)m$, even in very special cases; but we are able to provide such a formula when $M = S$ and 
$I = J+(x_0,\dots x_n)^D$ for some $D$, with $J$ generated in a single degree, in terms of the numbers $e_m(J,M)$. In particular, we find that in this situation the numbers
$e_m(I,M) - e_{m+1}(I,M)$ need not be weakly decreasing.

Section 3 of the paper uses some of the same ideas to prove a result close in spirit to the Eisenbud-Goto conjecture. Let $I$ 
be a reduced, equidimensional homogeneous ideal in $S$, and suppose that $k$ is algebraically closed. The Eisenbud-Goto conjecture then asserts the following:
\emph{ if the projective variety $X$ associated to $I$ is connected in codimension 1, then $\reg I \leq \deg X - \codim X +1$.} This conjecture 
is wide open, even for smooth varieties $X$, when the dimension of $X$ is large.

In the conjecture the hypothesis ``connected in codimension 1'' is necessary, as an example of Giaimo (included in Section 3) shows---
without the hypothesis, one must expect exponentially large regularity in general. But we are able to 
 prove a bound that is only slightly weaker than that of the Eisenbud-Goto conjecture \emph{without any connectedness hypothesis,} 
assuming instead that $X$ lies on a quadric (and has only isolated  ``bad'' singularities).

\section{$\gm$-Primary Ideals Generated in One Degree}
 
In this section, $S$ denotes a standard graded algebra over a field $k$.
We write $\gm$ for the homogeneous maximal ideal of $S$.
Let $I\subset S$
be a homogeneous ideal generated in a single degree $d$.

We consider the Rees ring $\cR (I) = S[It]$ of $I$, a standard
bigraded $k$-algebra as described above. 
Let $A$ be the ring
$$
A:=k[I_dt] = \oplus_j\cR(I)_{(0,j)}\subset \cR(I).
$$
It is a bigraded subalgebra of $\cR (I)$,
generated in degree $(0,1)$, which is a direct summand as an $A$-module.
We regard $A$ as a standard graded algebra, generated in degree 1 over $k$.
We write $\gn$ for the homogeneous
maximal ideal of $A$. Since $I$ is generated in one degree, $A$ is isomorphic
to the special fiber ring $\cF(I)=\cR(I) \otimes_S k$.

For $M$ a finitely generated graded $S$-module we consider 
the Rees module $\cR (I,M)=S[It]M$, which is a finitely generated bigraded $\cR(I)$-module.
We define 
$$
N_i(I,M):=k[I_dt]M_i\subset \cR(I,M).
$$ 
With the $(0,1)$-grading,  $N_i(I,M)$ is generated in degree 0, and has degrees
determined by the powers of $t$.
As an $A$-module, $\cR(I,M)$ is isomorphic to the direct sum of the $N_i(I,M)$.
In particular 
$$
\reg_{(y_0,\dots,y_r)} \cR(I,M),
$$
the $(0,1)$-regularity of $\cR(I,M)$, is the maximum
of the regularities of the $N_i(I,M)$ (as $A$-modules).
We shall see later how to restrict the range of $i$ required.

\begin{theorem}
\label{main}
Suppose that $I\subset S$ is an  ideal generated by forms
of a single degree $d$, and $M$ is a finitely generated graded $S$-module,
generated in a single degree,
such that $M/IM$ has finite length. Let $e$ be the number such that
$$
\reg I^mM=md+e
$$
for $m\gg 0$. Let $N_e=N_e(I,M)$.
\begin{enumerate}
\item The equality
$
\reg I^mM=md+e
$
holds if
$$
m\geq \max\{\reg H_\gn^1(N_e)+1,\   \frac{\reg M-e+1}{d}\}. 
$$
\item In case 
$
\reg H_\gn^1(N_e)\geq  {(\reg M-e+1)}/{d},
$
and $m\geq 1$, the equality 
$
\reg I^mM=md+e
$
holds if and only if 
$$
m\geq \reg H_\gn^1(N_e)+1.$$
\end{enumerate}
\end{theorem}

\begin{corollary}\label{simple bound}
Let $I,S,M,d,e$ be as in Theorem \ref{main}.
The equality
$\reg I^mM=md+e$
holds for all $m\geq \max\{\reg_y \cR(I,M),  \frac{\reg M+1}{d}\}. $
\end{corollary}

\begin{proof}[Proof of the Corollary] Since $N_e$ is an $A$-direct summand of $\cR(I,M)$,
$$
\reg H^1_\gn(N_e)+1 \leq \reg N_e \leq \reg_{(y_0,\dots,y_r)}\cR(I,M).
$$
\end{proof}

\begin{proof}[Proof of the Theorem]
Consider first part 1, and assume that 
$$
m\geq \max\{\reg H_\gn^1(N_e)+1,\  \frac{\reg M-e+1}{d}\}.
$$
By Eisenbud-Harris [2008] Proposition 1.1, $\{e_n\}$ is a non-increasing sequence. Thus it suffices to show that  
$\reg I^mM \leq md+e$. Our assumption on $m$ implies 
that $\reg M\leq md+e-1$.
Because of the exact sequence
\begin{equation}
0\to I^mM \to M\to M/I^mM\to 0 
\end{equation}
we only need to show that 
$\reg M/I^mM\leq md+e-1$. Since $M/I^mM$ has finite length,
this is equivalent to the statement that
$$
(I^mM)_{md+e}= M_{md+e}.
$$
The definition of $e$ implies, by the same argument,
that this equality at least holds for sufficiently large $m$.

Let $N'_e = N_e(\gm^d,M) = \oplus_{j\in \ZZ} M_{jd+e}t^j$. Note that
$N_e'$ is naturally a graded $A$-module (with $j$-th graded piece
$M_{jd+e}t^j$) and that $N_e$ is a submodule.
Let 
$$
E=N_e'/N_e=\frac{\oplus_jM_{jd+e}t^j}{\oplus_j(I^j)_{jd}M_et^j}.
$$
By the preceding remark, the module $E$
has finite length. 

We wish
to show that $E_m=0$.
Since $m\geq \reg H^1_\gn(N_e)+1$ we see from the exact sequence
\begin{equation}
\cdots\to H^0_\gn (N_e') \to E \to H^1_\gn(N_e)\to H^1_\gn(N_e')\to \cdots
\end{equation}
that it suffices to prove $H^0_\gn(N_e')_m=0$.

We may identify the $A$-module $N_e'$ with 
the $k[I_d]$-module $\oplus M_{dj+e}$, which is a $k[I_d]$-direct summand of $M$.
Note that this identification sends the degree $j$ part of $N_e'$ to the degree $dj+e$ 
part of $M$.
Moreover, since
$I_dS=I$ contains a power of $\gm$,
the module  $H^0_\gn(N_e')$ is a summand 
of $H^0_\gm(M)$ (with the same degree shift).
On the other hand, $H^0_\gm (M)_{dj+e}=0$ when $dj+e\geq 1+\reg M$.
Thus $H^0_\gn(N_e')_j=0$ when $j\geq (\reg M-e+1)/d$, concluding the proof of
part 1.

We now consider part 2. Given part 1 and Proposition \ref{decreasing}, it suffices to
show that if $m=\reg H^1_\gn (N_e)$ then $\reg I^mM\geq md+e+1$.
It follows from the hypothesis of part 2 that $\reg M\leq md+e-1.$
Because of the exact sequence (1) we only need to show that 
$\reg(M/I^mM)\geq md+e.$ Let $N_e'$ and $E$ be as in part 1. 
We want to show that $E_m\neq 0$.

Using exact sequence (2) and the fact that $H^1_\gn(N_e)_m\neq 0$, we see that it suffices to show 
$H^1_\gn(N_e')_m= 0.$ Since $N_e'$ is a summand of $M$ (with a shift of degree)
it suffices to show $H^1_\gm(M)_{md+e}=0$. This holds because, by hypothesis,
$\reg M\leq md+e-1$.
\end{proof}

\begin{conjecture} 
\label{decreasing conjecture}
If $I,S, M$ are as in Theorem \ref{main}, then the regularity of $N_i$ is non-increasing
from $i=0$. In particular, the $(0,1)$-regularity of $\cR(I)$ is equal to the regularity
of $k[I_d]$.
\end{conjecture}

We can prove the conjecture in the case where $I$ is a power of the maximal
ideal.

\begin{proposition} \label{yregN(md)}
Let $M$ be
a finitely generated graded $S$-module, generated in degree 0. 
\begin{enumerate}
\item If $i\geq 0$, then
$$
\reg N_i(\gm^d,M) \leq \max\biggl\{0, \frac{\reg M - i + (d-1)\dim M}{d} \biggr\}.
$$ 
In particular
$\reg N_i(\gm^d,M)=0$ for $i\geq \reg M+ (d-1)(\dim M - 1)$. 
\item If $H^0_\gm(M)=0$, then the inequality of part 1) is an equality. In particular,
 the sequence of numbers $\{\reg N_i(\gm^d,M) \mid i\geq 0\}$
is weakly decreasing.
\end{enumerate}
\end{proposition}

\begin{proof}
In the previous proof we have seen that there is a homogeneous isomorphism of $k[S_d]$-modules 
$$
N_i \cong M_iK[S_d](i)=\oplus_{j\geq0}M_{dj+i}=(M(i)_{\geq 0})^{(d)},
$$
where we consider $N_i$ as a $k[S_d]$-module via the identification 
$k[S_dt] \cong k[S_d]$; here $-^{(d)}$ denotes the Veronese functor.

The exact sequence
$$
0\to M(i)_{\geq 0} \to M(i) \to M(i)/M(i)_{\geq 0}\to 0
$$
gives rise to an exact sequence
\begin{align*}
0\to  &H^0_\gm (M(i)_{\geq 0}) \to H^0_\gm(M(i))\to M(i)/M(i)_{\geq 0}\\
&\to  H^1_\gm (M(i)_{\geq 0}) \to H^1_\gm(M(i))\to 0
\end{align*}
and isomorphisms $H^\ell_\gm (M(i)_{\geq 0}) \cong H^\ell_\gm(M(i))$
for $2\leq \ell$.

Since the $d$-th Veronese functor commutes with taking local cohomology it follows that 
\begin{align}
&\reg (M(i)_{\geq 0})^{(d)}\notag \\
&\leq\max\{1+\reg (M(i)/M(i)_{\geq 0})^{(d)}, \max\{ \reg (H_\gm^\ell (M(i)))^{(d)}+\ell \mid 0\leq \ell\leq \dim M\}\notag\\
&=\max\biggl\{0, \max\biggl\{\biggl\lfloor \frac{\reg H^\ell (M)-i}{d}\biggr\rfloor+\ell \mid 0\leq \ell\leq \dim M\biggr\}\biggr\}\notag\\
&\leq \max\biggl\{0, \max\biggl\{\biggl\lfloor  \frac{\reg M-i-\ell}{d}\biggr\rfloor +\ell \mid 0\leq \ell\leq \dim M\biggr\}\biggr\}\notag\\
 & \leq\max\biggl\{0,  \biggl\lfloor \frac{\reg M-i+(d-1)\dim M}{d} \biggr\rfloor\biggr\}
\end{align}
which gives the desired formula. If $H^0_\gm(M) = 0$ then
 the first inequality is an equality, which implies part 2.
\end{proof}

We can also prove Conjecture \ref{decreasing conjecture} for $i\geq e$, at least when
$H^0_\gm(M)=0$.
\begin{proposition}
\label{decreasing}
Suppose that $I\subset S$ is an  ideal generated by forms
of a single degree $d$, and $M$ is a finitely generated graded $S$-module,
generated in a single degree,
such that $M/IM$ has finite length and $H^0_\gm(M)=0$. For each $m$, let $e_m$ be the number such
that $\reg I^mM = md+e_m$, and let $e=e_m$ for $m\gg 0$. Let $N_j$ be the module
defined above.
\begin{enumerate}
\item $e_m\geq e_{m+1}\geq e_m-d$.
\item If 
$i\geq e$ then $\reg N_{i+1}\leq \reg N_i$.
\end{enumerate}

\end{proposition}

\begin{proof}
The inequality $e_m\geq e_{m+1}$
of part 1 is proven in Eisenbud-Harris \cite{EH}, Proposition 1.1.

For the second inequality it suffices to prove that 
$\reg I^mM \leq \reg I^{m+1}M$, for then
$dm +e_m \leq d(m+1) +e_{m+1}$, that is,
$e_m \leq d+e_{m+1}$.

Recall that $M/I^{m+1}M$ has finite length and
$H^0_{\gm}(M)=0$. The exact sequence
$$
0\to I^{m+1}M \to M\to M/I^{m+1}M \to 0
$$
shows that 
$ \reg H^1_{\gm}(I^{m+1}M) = \max \{ \reg M/I^{m+1}M , \reg H^1_{\gm}(M) \}$
and moreover  $H^\ell_{\gm}(I^{m+1}M) = \reg H^\ell_{\gm}(M)$ for $\ell\geq 2$.
The same equalities hold for $I^mM$ in place of $I^{m+1}M$. 
The epimorphism of
finite length modules $M/I^{m+1}M \twoheadrightarrow M/I^{m}M$
implies that $\reg M/I^{m+1}M \geq \reg M/I^{m}M$, 
and the desired inequality follows.

For part 2, we note that for $i\geq e$ we can embed
$N_i$ into $N'_i:=N_i(\gm^d,M)$ with finite length cokernel.
From $H^0_\gm(M)=0$ we deduce $H^0_\gm(N_i')=0$ and thus $H^0_\gm(N_i)=0$.
Therefore
$
\reg N_i = \max\{ \reg N_i', \reg (N_i'/N_i)+1\}.
$

Since $H^0_\gm(M)=0$, part 2 of Proposition \ref{yregN(md)} shows that the 
numbers $\reg N_i'$ are weakly decreasing.
On the other hand, the generators of $\gm$
provide an epimorphism $\oplus N_i' \to N_{i+1}'$
that induces an epimorphism
$\oplus N_i'/N_i \to N_{i+1}'/N_{i+1}'$.  Thus the $(0,1)$ regularity of
the finite length module $(N_i'/N_{i+1})$ is also weakly decreasing when $i\geq e$.
\end{proof}


\begin{corollary}
Let $S=k[x_1,\dots, x_n]$
and let $I,d,e$ be as in Theorem \ref{main}.
If $e=0$ and $m\geq \reg k[I_d]$, then 
$\reg I^m=md+e$.
\end{corollary}
\begin{proof}
One uses Theorem \ref{main}.1 and Proposition \ref{decreasing}.2.
\end{proof}



\begin{example}
The regularity of $\cR(I)$ is often much larger than
the regularity of the module $N_e$.
For the ideal 
$I=(x^{20}, x^3y^{17}, x^{12}y^8, y^{20})\subset k[x,y]$
we have $\reg I^m \geq 20m+7$, with equality if and only if $m\geq 2$.
Here the $(0,1)$ regularity of the Rees algebra, and also the regularity of
$k[I_d]$, are equal to 7. By Theorem \ref{main}, 
$\reg H^1_\gn (N_e) \leq 1$ (and in fact equality holds).

For the ideal 
$I=(x^{20}, x^3y^{17}, x^{25}y^5, y^{20})\subset k[x,y]$
we have $\reg I^m \geq 20m+4$, with equality if and only if $m\geq 4$.
Here again the $(0,1)$ regularity of the Rees algebra, and also the regularity of
$k[I_d]$, are equal to 7. By Theorem \ref{main}, 
$\reg H^1_\gn (N_e) \leq 3$ (and again, in fact, equality holds).

\end{example}

\section{Ideals With Generators in More Than One Degree}

As a first example, we have:

\begin{proposition} 
Let $I\subset S=k[x_1,\dots,x_n]$ be a homogeneous ideal,
and $M$ a finitely generated graded $S$-module.
If $I\subset S$ is generated by an $M$-regular sequence
of degrees $d=d_1\geq \cdots\geq d_t$ and $m\geq 1$ then
$\reg I^mM = dm+e$ where
$
e=\reg M+ \sum_{i=2}^t (d_i-1).
$
\end{proposition}

\begin{proof} Since $I$ is generated by a regular sequence on $M$, we may tensor $M$ with the Eagon-Northcott resolution
of $I^m$ and get a resolution of $I^m\tensor M = I^mM$ by shifted copies of $M$.  
Analyzing the shifts, we see that $\reg I^mM = dm+e$.
\end{proof}

\begin{corollary} Let $I\subset S=k[x_1,\dots,x_n]$ be a homogeneous ideal,
and $M$ a finitely generated graded $S$-module.
Let $d$ be the asymptotic generator degree of $I$ on $M$, and write 
$\reg I^mM = dm+e_m$. 
If $I$ contains an $M$-regular sequence
of degrees $d=d_1\geq \cdots\geq d_t$ with $t=\dim M$,
then
$
e_m\leq \reg M+ \sum_{i=2}^n (d_i-1).
$
for every $m \geq 1.$
\end{corollary}

\smallskip

In general, we can analyze only special cases. 

\begin{theorem}\label{mixed degrees}
Let $J\subset S=k[x_1,\dots,x_n]$ be an $\gm$-primary ideal generated
by forms of a single degree $d$. Write $I=J+\gm^{d+k}$ for some $k\geq 0$.
Let $f_m(p) = (d+k)m-kp$, and
$$
p_m = \min\{p\geq 1 \mid \reg J^p \geq f_m(p)\}.
$$
For $m\geq 1$ we have
$$
\reg I^m = \min\{\reg J^{p_m}, f_m(p_m-1)\}.
$$

\end{theorem}

\begin{proof}
Define $e_p$ by the formula
$\reg J^p = dp+e_p$. Note that $p_m$ is finite, and in fact
$p_m \leq m$ since $\reg J^m \geq dm$.

We have 
$$
I^m = \sum_{p=0}^m J^p(\gm^{d+k})^{m-p}.
$$
Thus, $\reg I^m \leq \min \{\reg J^p(\gm^{d+k})^{m-p}\mid 0\leq p\leq m\}.$
Moreover,
$J^p(\gm^{d+k})^{m-p}= (J^p)_{\geq dp+(d+k)(m-p)} = (J^p)_{\geq f_m(p)}$, so
$$
\reg J^p(\gm^{d+k})^{m-p}=\max\{ \reg J^p, f_m(p)\}.
$$

We claim that the minimum value of $\reg J^p(\gm^{d+k})^{m-p}$ is taken on
either for $p=p_m$ or $p=p_m-1$, and that in either case it is
$$
\min_{0\leq p\leq m} \{\reg J^p(\gm^{d+k})^{m-p}\} 
= \min\{ \reg J^{p_m}, f_m(p_m-1)\}.
$$
This follows because, as $p$ increases,
 the function $\reg J^p$ is weakly increasing while
$f_m(p)$ is decreasing, and
for $p=m$ the first is at least as large as the second, and $p_m\geq 1$---see Figure 1.
Note that the minimum value is the value claimed in the Theorem
for $\reg I^m$. 
\begin{figure}
\centering\includegraphics[scale=.5]{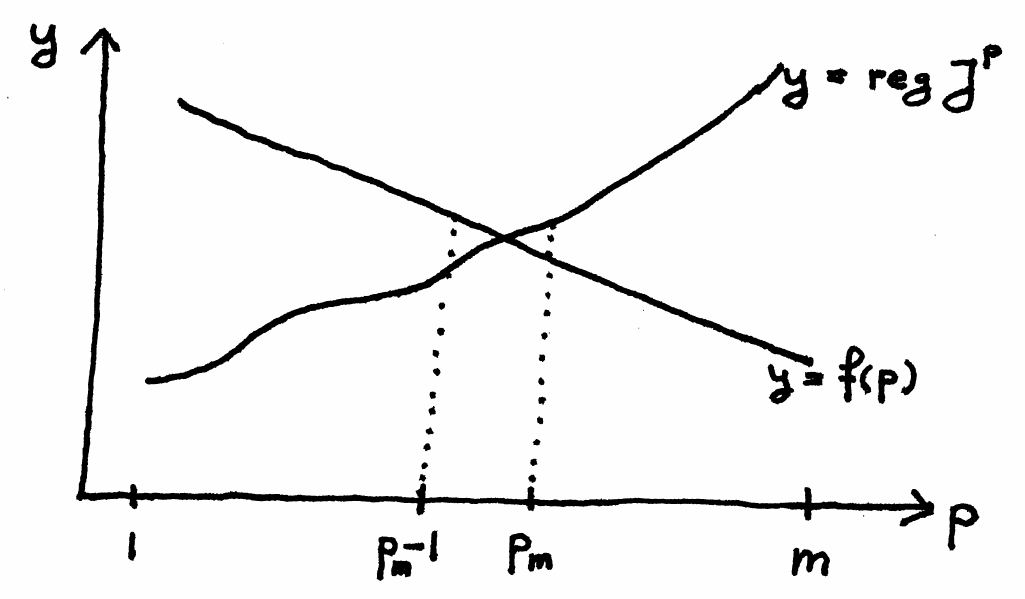}
\caption{Where the graphs of $f_m(p)$ and $\reg J^p$ cross}
\end{figure}
\goodbreak

Thus it is enough to show that 
$$
\reg I^m \geq  \min\{ \reg J^{p_m}, f_m(p_m-1)\}.
$$
Write $a=\min\{ \reg J^{p_m}, f_m(p_m-1)\}.$
Note that $I^m\subset J^{p_m}+\gm^{f_m(p_m-1)}$. 
Thus it suffices to prove that 
$$
\gm^{a-1}\not\subset J^{p_m}+\gm^{f_m(p_m-1)}.
$$
Since $a-1< f_m(p_m-1)$, this is equivalent to $\gm^{a-1} \not \subset J^{p_m}$.
But the latter holds because $a-1< \reg J^{p_m}.$
\end{proof}

\begin{example} If $I$ is not generated in a single degree then in the formula 
$\reg I^m = md +e_m$ the $e_m$ may not be weakly
decreasing. It can even go up and then down. For example, 
using Theorem \ref{mixed degrees} one 
can easily compute that
 if
$$
I=(x_1^4,\dots,x_4^4)(x_1,\dots,x_4)+(x_1,\dots,x_4)^6 \subset S=k[x_1,\dots,x_4]
$$
then $\reg I^m = 5m+e_m$, where the successive values of $e_m$  for $m= 1,2,\dots$
are
$1,2,2,1, 1,1,1,1,0,0,0,\dots$.
\end{example}

\begin{proposition}\label{depth1 graded ring}
Let $I\subset S=k[x_1,\dots,x_n]$ be a homogeneous ideal,
and $M$ a finitely generated graded $S$-module, concentrated in non-negative
degrees, such that $M/IM$ has finite length.
Let $d$ be the asymptotic generator degree of $I$ on $M$, and write 
$\reg I^mM = dm+e_m$. 
\begin{enumerate}
\item If $I$ is generated in degrees $\leq d$, then the sequence of integers
$\{e_m\mid m\geq (\reg M+1)/d\}$ is weakly decreasing.
\item If the associated graded module
$\gr_I(M)$ has positive depth, then the sequence $\{e_m\mid m\geq \reg M /d\}$
is weakly increasing. 
\end{enumerate}
\end{proposition}

\begin{proof} We first prove part 1.
 If $I$ is generated by homogeneous elements of degrees $d_i$
then multiplication by these elements gives a surjection
$$
\oplus_i \bigg(\frac{I^{m-1}M}{I^mM}(-d_i)\bigg) \to \frac{I^mM}{I^{m+1}M}
$$
of modules of finite length. Thus 
$$
\reg I^mM/I^{m+1}M\leq \reg I^{m-1}M/I^mM + d\leq \reg M/I^mM+d.
$$
Now the exact sequence 
$$
0\to I^mM/I^{m+1}M \to M/I^{m+1}M \to M/I^mM \to 0
$$
shows that $\reg M/I^{m+1}M\leq \reg M/I^mM +d.$

Since $\reg (I^m)^p M = (dm)p + e_{mp}$ for $p \gg 0$, we conclude that
the asymptotic generator degree of $I^m$ on $M$ is $dm$. Thus the
generator degree of $I^mM$ is at least $dm$ because $M$ is
concentrated in non-negative degrees. It follows that $\reg I^mM \geq dm$.
Thus, 
if $m\geq (\reg M+1)/d$ then $\reg M \leq dm-1\leq \reg I^mM-1$. Now
the inequality $\reg M/I^{m+1}M\leq \reg M/I^mM +d$
implies that $\reg I^{m+1} M \leq \reg I^mM+d$.

\smallbreak
For part 2
we may assume that $k$ is infinite. The definition of $d$ shows that
for some integer $p$ we have
$$
(I/I^2)^p \gr_I(M) \subset ((I_{\leq d} + I^2)/I^2)\gr_I(M).
$$
It follows that there exists an element $a\in I_d$
whose leading form $a+I^2\in \gr_I(S)$ is a non-zerodivisor on $\gr_I(M)$.
Hence $I^{m+1}M:_Ma = I^mM$. Thus multiplication by
$a$ induces an embedding
$$
\frac{M}{I^mM}(-d)\hookrightarrow \frac{M}{I^{m+1}M}. 
$$
On the other hand, 
$a(I^mM:_M\gm)\subset  I^{m+1}M:_M\gm$. 
Therefore
$$
\frac{I^mM:_M\gm}{I^mM}(-d) \hookrightarrow \frac{I^{m+1}M:_M\gm}{I^{m+1}M}.
$$
This implies $\reg M/I^{m+1}M\geq \reg M/I^mM+d$, and hence $\reg I^{m+1}M\geq \reg I^mM+d$
whenever $m\geq \reg M /d$. 
\end{proof}

\begin{corollary}\label{constancy}
 Let $I\subset S=k[x_1,\dots,x_n]$ be a homogeneous $\gm$-primary ideal
with asymptotic generator degree $d$. If $I$ is generated in degrees $\leq d$ and $\gr_I(S)$
has positive depth, then $\reg I^m = dm+e$ for some $e$ and every $m\geq 1$. \qed
\end{corollary}

\begin{example}
One cannot drop the assumption of generation in degree $\leq d$ from Corollary \ref{constancy}.
If
$$
I=(x^4,y^4,z^4)+(x,y,z)^5 \subset S=k[x,y,z],
$$
then $\reg I^m = 4m+e_m$, where the successive values of $e_m$  for $m= 1,2,\dots$
are
$1,2,2,2,2,\dots$.
Computation with Macaulay2 shows that the depth of the associated
graded ring of $I$ is at least 1. 
\end{example}


\section{A Case of the (Almost) Eisenbud-Goto Conjecture}

Eisenbud and Goto [1984] conjecture that the regularity of
a nondegenerate, geometrically
 reduced irreducible subscheme $X\subset \PP^n$
 has regularity at most $\deg X-\codim X+1$.
 They further conjecture that the hypothesis can be weakened to
 say that the nondegenerate scheme is geometrically reduced and connected in codimension 1,
 and this has been proved by Giaimo [2006] for curves.
 The bound can fail for disconnected schemes. For example, if $X$ is the union
 of two skew lines in $\PP^3$ then the degree of $X$ is 2 but the regularity
 (that is, the regularity of the ideal of $X$) is 2 rather than 1. Derksen and Sidman [2002] 
 have shown that in general a union of linear subspaces of projective space has regularity
 at  most the number of subspaces. 
 
 One might guess from this that the regularity of a reduced equidimensional
 scheme would be bounded by the degree of the scheme, but this is not the case.

\begin{example} Here is a reduced equidimensional union of two irreducible complete intersections whose regularity is
much larger than its degree: 

By Mayr-Meyer [1982] there is a homogeneous ideal $I\subset S=\CC[x_1,\dots,x_n]$ generated by
$10n$ forms of degrees two and three, having regularity of the order of $2^{2^n}$. In the ring
$R=S[z_1,\dots]$ we build an ideal $I'$ whose generators correspond to those of $I$
by replacing the monomials in the generators of $I$ with products of new variables $z_j$ in such
a way each $z_j$ occurs only linearly, and no $z_j$ occurs twice. Clearly the generators of this new
ideal are a regular sequence. If any of the generators are monomials, we add further new variables $w_j$
and make each a binomial that will be a prime. Since the variables are all distinct, the resulting
complete intersection will also be prime, and modulo an ideal of the form
$L=(\{z_j - x_{p(j)}\} + (\{ w_j\})$ the ideal $I'$ becomes equal to the ideal $I$. The codimension
of $L$ is clearly at least as big as the codimension of the complete intersection. We add further linear
forms to the complete intersection $I'$ to make the codimensions the same. 

The ideal $I'\cap L$ now defines the union of two reduced, irreducible complete intersections, while
the ideal $I'+L$ defines the original Mayr-Meyer example. From the short exact sequence 
$$
0\to I'\cap L \to I'\oplus L \to I'+L\to 0
$$
we see that the regularity of $I'\cap L$ is of the order of $2^{2^n}$. On the other hand, the degree of the
subscheme defined by $I'\cap L$ is at most of the order of $3^{10n}$.
\end{example} 
 
  We state our result in terms
 of the regularity of the homogeneous coordinate ring $S$ of $X$, which is 
 one less than $\reg X$,  to emphasize the parallel between the two parts 
 of the Theorem. Recall that a local algebra essentially of finite type over a field characteristic zero is said to have a rational
 singularity if it is normal and Cohen-Macaulay and, if $\pi: \tilde X \to \Spec R$ is a resolution of singularities, 
 then $\pi_*(\omega_{\tilde X}) = \omega_{\Spec R}$.

\begin{theorem}
Let $X$ be a reduced equidimensional subscheme of codimension 2 in $\PP^n_k$ where $k$ is a field
of characteristic zero and the locus of non-rational singularities of $X$
has dimension zero. Let $S_X$ be the homogeneous coordinate ring of $X$.
If $X$ lies on a quadric hypersurface, then 
\begin{enumerate}
\item
$\reg S_X\leq \deg X$.  
\item 
If $x_1,\dots,x_n$ are general linear forms in $S_X$,
 and $I$
is the ideal they generate, then $\reg_y\cR(I,S_X)\leq \deg X-\codim X +1$.
\end{enumerate}
\end{theorem}

Note that the Eisenbud-Goto conjecture would say, under the additional hypothesis that $X$ is connected in codimension 1, that
$\reg S_X \leq \deg X-\codim X  = \deg X - 2$.

 \begin{proof} We make use of the notation introduced in part 2 of the Theorem, and
we write $\gm$ for the homogeneous maximal ideal of $S_X$.
 Let $\cF= k[I_1]\subset S_X$ and note that  $\cF$ is isomorphic to the fiber ring
 $\cF \cong \cR(I,S_X)/\gm \cR(I,S_X)$.
 Let $x$ be a linear form such that $\gm  = (I,x)$.
 Because the $x_1,\dots,x_n$ are general and the ideal defining $X$ contains a quadric,
 $S_X=\cF+\cF x$. Thus $S_X/\cF\cong  (\cF/(\cF:_\cF S_X))(-1)$. The extension
$\cF \subset S_X$ is birational. Hence $\cF$ is the 
 ring of a hypersurface whose degree is $ \deg S_X$ in $\PP^{n-1}$. It follows that $\reg \cF = \deg S_X -1$.
 
As $\omega_\cF = \cF(-n+\deg S_X)$ we have
$\cF:_\cF S_X = \Hom_\cF(S_X,\cF) = \omega_{S_X}(n-\deg S_X)$. The hypothesis that the 
characteristic is zero and that the equidimensional scheme  $X$
has at most isolated non-rational singularities
implies that the regularity of $\omega_{S_X}$ is at most $\dim S_X = n-1$ (see Chardin-Ulrich
\cite{Chardin-Ulrich} Theorem 1.3,
which is based on results of Ohsawa [1984] and Koll\'ar [1986], Theorem 2.1(iii)). 
It follows that
$\reg (\cF:_\cF S_X) \leq n-1 - (n -\deg S_X) = \deg S_X -1.$ Thus 
$\reg S_X/\cF \leq \deg S_X$, and therefore $\reg S_X \leq \deg S_X$, proving the first statement.

For the second statement, let $G=\gr_IS_X$ be the associated graded ring of $S_X$ with respect to $I$, which 
is an $S_X$-module via the  map $S_X\to S_X/I = G_0$. 
By 
Johnson and Ulrich [1996] Proposition 4.1 one has
$\reg_y\cR(I,S_X) = \reg_y  G$, so it suffices to bound the latter.

Note that $\cF = G/\gm G = G/xG$. Because the ideal defining $X$ contains a quadric we have $x^2\in I$.
It follows that $x^2 G = 0$. Of course $xG\cong G/(0:_G x)$. We will show that 
$G/(0:_G x)\cong \cF/(\cF:_\cF S_X)$. Indeed, the inclusion $\cF\subset \cR(I,S_X)$
induces an inclusion $\cF\subset G$, and hence a map $\cF \to G/(0:_G x)$ which
is surjective because $xG\subset 0:_G x.$ To compute the kernel, let $f\in \cF$ be
a form of degree $i$. We have $fx = 0 $ in $G$ if and only if, as elements of $S_X$, 
we have $fx \in I^{i+1}$. But the degree (in $S_X$) of $fx$ is $i+1$, so 
this happens if and only if $fx \in \cF_{i+1}$. This in turn means that 
$f\in \cF:_\cF x = \cF:_\cF S_X$.

From the computation of the regularity of $\cF:_\cF S_X$ above, we get
$\reg G \leq \max\{ \reg \cF/(\cF:_\cF S_X), \reg \cF\} = \deg S_X -1$.
\end{proof}

\vfill\eject

\bigskip

\vbox{\noindent Author Addresses:\par
\smallskip
\noindent{David Eisenbud}\par
\noindent{Department of Mathematics, University of California, Berkeley,
Berkeley, CA 94720}\par
\noindent{eisenbud@math.berkeley.edu}\par
\smallskip
\noindent{Bernd Ulrich}\par
\noindent{Department of Mathematics, Purdue University, West Lafayette, IN 47907}\par
\noindent{ulrich@math.purdue.edu}\par
}

\end{document}

%% file: preamble.tex
\def\antiddot{\mathinner{\mkern1mu\raise1pt\vbox{\kern7pt\hbox{.}}\mkern2mu
        \raise4pt\hbox{.}\mkern2mu\raise7pt\hbox{.}\mkern1mu}}

\newcommand{\CC}{{\mathbb C}}

\newcommand{\PP}{{\mathbb P}}

\newcommand{\ZZ}{{\mathbb Z}}

\newcommand{\ann}{{\rm{ann}}}

\newcommand{\s}{\mathcal}

\newcommand{\cF}{{\s F}}

\newcommand{\cR}{{\s R}}

\newcommand{\tensor}{\otimes}

\newcommand{\punkt}{\hspace{-.3ex}\raise.15ex\hbox to1ex{\Huge.}}

\def \fix#1 {{\hfill\break \bf (( #1 ))\hfill\break}}

\DeclareMathOperator{\reg}{reg}
\DeclareMathOperator{\Spec}{Spec}

\DeclareMathOperator{\Hom}{Hom}

\DeclareMathOperator{\codim}{codim}

\newcommand{\gm}{\mathfrak m}
\newcommand{\gn}{\mathfrak n}
\def\gr{{\mathfrak {gr}}}

\newtheorem{theorem}{Theorem}[section]

\newtheorem{proposition}[theorem]{Proposition}
\newtheorem{corollary}[theorem]{Corollary}
\newtheorem{conjecture}[theorem]{Conjecture}
\theoremstyle{definition}

\newtheorem{example}[theorem]{Example}